\documentclass[11pt,reqno,twoside]{article}
\usepackage{amsmath,amssymb,theorem,color,epsfig,epic,subfigure,amsfonts,hhline,graphicx}
\usepackage[normalem]{ulem}
\usepackage{cancel}
\usepackage[colorinlistoftodos]{todonotes}
\theorembodyfont{\itshape}
\newtheorem{thm}{Theorem}[section]
\newtheorem{lem}[thm]{Lemma}
\newtheorem{prop}[thm]{Proposition}
{\theorembodyfont{\upshape}
\newtheorem{define}[thm]{Definition}
\newtheorem{rem}[thm]{Remark}
\newtheorem{ex}[thm]{Example}
}
\newtheorem{cor}[thm]{Corollary}

\newcommand{\Proof}[1][]{\noindent{\itshape Proof#1. }}
\newcommand{\EndProof}{\hfill$\Box$\bigskip}

\makeatletter
\let\@fnsymbol\@arabic
    
    \@addtoreset{equation}{section}
\makeatother

\def\hat{\widehat}
\def\tilde{\widetilde}

\def\sbs{\subset}

\def\a{\alpha}

\newcommand\unit{\hbox{\rm 1\kern-2.8truept l}}

\newcommand{\tr}[1]{\hbox{\rm{tr}} (#1)}

\newcommand\Lform{{\mathcal{L}}\kern-7.56pt\raise1.0pt\hbox{$-$}}

\begin{document}
\title{Group of automorphisms for strongly quasi invariant states}
\author{ Ameur Dhahri\footnote{Dipartimento di Matematica, Politecnico di Milano, Piazza Leonardo di Vinci 32, I-20133 Milano, Italy. E-mail: ameur.dhahri@polimi.it}, Chul Ki Ko\footnote{University College, Yonsei University, 85 Songdogwahak-ro,
Yeonsu-gu, Incheon 21983, Korea. E-mail: kochulki@yonsei.ac.kr}, and Hyun Jae Yoo\footnote{Department of Applied Mathematics and Institute for Integrated Mathematical Sciences,
Hankyong National University, 327 Jungang-ro, Anseong-si,
Gyeonggi-do 17579, Korea. E-mail: yoohj@hknu.ac.kr}}
\date{ }
   \maketitle

\begin{abstract}
For a $*$-automorphism group $G$ on a von Neumann algebra, we study the $G$-quasi invariant states and their properties. The $G$-quasi invariance or $G$-strongly quasi invariance are weaker than the $G$-invariance and have wide applications. We develop several properties for $G$-strongly quasi invariant states. Many of them are the extensions of the already developed theories for $G$-invariant states. Among others, we consider the relationship between the group $G$ and modular automorphism group, invariant subalgebras, ergodicity, modular theory, and abelian subalgebras. We provide with some examples to support the results.
\end{abstract}

\section{Introduction}
In this paper we discuss the group of automorphisms and quasi invariant states thereof. The closed quantum dynamical system is often represented with a group of $*$-automorphisms on a $C^*$- or von Neumann algebra $\mathcal A$. The invariant states, or the equilibrium states for the group play the central role in the understanding of the system. From the late sixties and early seventies many theories have been developed and found applications in quantum statistical mechanics. See \cite{BR, DKKR, He, HeTa, St, Ta} and references therein.

In addition to the time evolution, one can further consider the other $*$-automorphisms for spatial movements such as space translations or rotations. Given such a group $G$ of $*$-automorphisms, the natural question is to find the $G$-invariant states and develop the properties of them. The invariant states as well as the relationships between the group actions of $G$ and the modular automorphism group are studied in many papers, see for instance \cite{He, HeTa, St, Su, Ta}. For the automorphism group $G$, as a weaker condition for invariance, Accardi and Dhahri introduced a concept of quasi invariance \cite{AcDh}. It considers a kind of Radon-Nikodym derivatives for the states translated by the group action. Related concepts can also be found in some other papers \cite{Gu, Ko, PeTa}. In \cite{AcDh}, after giving the definition for $G$-quasi and $G$-strongly quasi invariance of states, the authors developed many concrete properties such as cocycle properties of the Radon-Nikodym derivatives. The purpose of this paper is to further develop the properties of $G$-quasi invariance.

Let us give a brief overview of this paper. In Section \ref{sec:general_theories}, we recall the definition of $G$-quasi invariance. We characterize the $G$-invariance of a state and the relationship between the modular automorphism group and $G$ (Theorem \ref{thm:G-invariance}). Then we study the invariant subalgebras (Proposition \ref{prop:mean_and_modular_group}) and ergodic properties (Theorem \ref{thm:ergodicity}). We consider several examples to support the theories developed. In Section \ref{sec:tracial}, we characterize a relationship between $G$-strongly quasi invariant states and tracial states, Theorem \ref{thm:characterization_by_tracial_states_amenable_group}. In Section \ref{sec:modular}, we develop the modular theory for the $G$-strongly quasi invariant states (Theorem \ref{thm:modular operator relation}). Finally in Section \ref{sec:abelian}, we discuss the abelian projections for the $G$-strongly quasi invariant states, Theorem \ref{thm:abelianness_compact_case}.

 \section{Automorphism group and quasi invariant states}\label{sec:general_theories}
 \subsection{$G$-invariance and $G$-strongly quasi invariance}
We briefly review some definitions and properties of quasi invariant states, respectively, of strongly quasi invariant states with respect to an automorphism group on a von Neumann algebra from the reference \cite{AcDh}. Let $\mathcal A$ be a von Neumann algebra and let $G$ be a topological group which has a strongly continuous representation $\{\tau_g:g\in G\}$ of normal $*$-automorphisms on $\mathcal A$. For notational simplicity, we write $g(a)$ for $\tau_g(a)$, $a\in \mathcal A$, and the strong continuity means that
\[
\lim_{g\to e} g(a)=a,\quad a\in \mathcal A,
\]
where $e\in G$ is the identity element of $G$.
\begin{define} \label{def:quasi_invariance}
A faithful normal state $\varphi$ on $\mathcal A$ is said to be $G$-quasi invariant if for all $g\in G$ there exists $x_g\in \mathcal A$, a Radon-Nikodym derivative, such that
\begin{equation}\label{eq:quasi_invariance}
\varphi(g(a))=\varphi(x_ga),\quad a\in \mathcal A.
\end{equation}
We say that $\varphi$ is $G$-strongly quasi invariant if it is $G$-quasi invariant and the Radon-Nikodym derivatives $x_g$ are self-adjoint: $x_g=x_g^*$ for all $g\in G$.
\end{define}
Sometimes it is convenient to introduce a generalized concept.
\begin{define}\label{eq:generalized}
We say that $\varphi$ is $G$-strongly quasi invariant in the generalized sense if for all $g\in G$ the relation \eqref{eq:quasi_invariance} holds for an $x_g=x_g^*$ affiliated with $\mathcal A$.
\end{define}
\begin{rem}\label{rem:properties} (i) When $\varphi$ is $G$-quasi invariant, the map $g\mapsto x_g$ is a normalized multiplicative left $G-1$-cocycle satisfying
\begin{equation}\label{eq:cocycle}
x_e=\unit,\quad x_{g_2g_1}=x_{g_1}g_1^{-1}(x_{g_2}), \quad g_1,g_2\in G,
\end{equation}
and $x_g$ is invertible with its inverse
\[
x_g^{-1}=g^{-1}(x_{g^{-1}}).
\]
Furthermore, it holds that $\varphi(x_ga)=\varphi(ax_g^*)$ for all $a\in \mathcal A$ and $g\in G$. \\
(ii) In the case that $\varphi$ is $G$-strongly quasi invariant in the generalized sense, the cocycle relation \eqref{eq:cocycle} still holds.\\
(iii) If $\varphi$ is $G$-strongly quasi invariant (or in the generalized sense), it holds that $x_g$ is strictly positive. In the case of $G$-strongly quasi invariance, It also holds that $x_g$ commutes with $x_h$ for all $g,h\in G$. Hence the $C^*$-algebra $\mathcal C$ generated by $\{x_g:g\in G\}$ is commutative. It can be also easily shown that $\mathcal C$ is a subalgebra of $\mathrm{Centr}(\varphi)$, the centralizer of $\varphi$, which is defined in \eqref{eq:centralizer} below. See \cite{AcDh} for the details.
\end{rem}
Recall that given a faithful normal state $\varphi$ on a von Neumann algebra $\mathcal A$, the set
\begin{equation}\label{eq:centralizer}
\mathrm{Centr}(\varphi):=\{x\in \mathcal A:\varphi(xy)=\varphi(yx)\text{ for all }y\in \mathcal A\}
\end{equation}
is called the centralizer of $\varphi$ \cite{He}.

Let $G$  be a  compact group with normalized Haar measure $dg$ and let $\varphi$ be a $G$-strongly quasi-invariant state on $\mathcal A$. We assume the map $g\mapsto x_g$ is continuous and define
\begin{equation}\label{eq:rnd_mean_and_Umegaki}
\kappa:=\int_G x_g dg\text{ and }E_G:=\int_Ggdg.
\end{equation}
Here the integral for $E_G$ is understood as an operator acting on $\mathcal A$, so $E_G(a)=\int_Gg(a)dg$ for all $a\in \mathcal A$.
Then, $\kappa$ is invertible with a bounded inverse \cite[Theorem 1]{AcDh}. The $*$-map $E_G$ is a Umegaki conditional expectation from $\mathcal A$ to $E_G(\mathcal A)=\mathcal F(G):=\{a\in \mathcal A:\,g(a)=a,\,\forall \,g\in G\}$, the $G$-invariant, or $G$-fixed subalgebra \cite[Theorem 2]{AcDh}. We define $\varphi_G(x):=\varphi(\kappa x)=\varphi(E_G(x))$, or equivalently $\varphi(x)=\varphi_G(\kappa^{-1}x)$. Particularly, using the translation invariance of the Haar measure, it is easily checked that $\varphi_G$ is a $G$-invariant state.
\begin{lem}\label{lem:centralizers1}
$\kappa, \kappa^{-1}\in \mathrm{Centr}(\varphi)\cap\mathrm{Centr}(\varphi_G)$.
\end{lem}
\Proof
Since $x_g\in \mathrm{Centr}(\varphi)$, it is obvious that $\kappa, \kappa^{-1}\in \mathrm{Centr}(\varphi)$. From this the fact $\kappa, \kappa^{-1}\in\mathrm{Centr}(\varphi_G)$ also easily follows.
\EndProof
\begin{rem}\label{rem:modular_group}
Let $(\sigma_t^\varphi)$ and $(\sigma_t^G)$ be the modular automorphism groups of $\varphi$ and $\varphi_G$, respectively. It holds that \cite[Theorem 4.6]{PeTa}
\[
 \sigma_t^\varphi(x)=\kappa^{-it}\sigma_t^G(x)\kappa^{it}, \quad t\ge 0,\,\,x\in \mathcal A.
\]
Note that $\varphi$ is invariant not only for $(\sigma_t^\varphi)$ but also for $(\sigma_t^G)$.
\end{rem}
Combining the results of \cite[Lemma 1]{HeTa} and \cite[Theorem 4.6]{PeTa} we can characterize the strong quasi invariance.
 \begin{thm}\label{thm:G-strong_quasi_invariance}
Suppose that the von Neumann algebra $\mathcal A$ is a factor.  A state $\varphi$ is $G$-strongly quasi invariant if and only if for all $g\in G$ there exists a strictly positive operator $x_g\in \mathcal A$ such that for all $a\in \mathcal A$,
  \begin{equation}\label{eq:G-strongly_quasi_invariance}
g^{-1}\circ\sigma_t^\varphi\circ g(a)=x_g^{it}\sigma_t^\varphi(a)x_g^{-it}.
 \end{equation}
 \end{thm}
 \Proof
 ($\Longrightarrow$) Suppose that $\varphi$ is $G$-strongly quasi invariant. We have $\varphi_g:=\varphi\circ g$ is a normal faithful state. By Lemma 1 of \cite{HeTa}, it holds that
 \begin{equation}\label{eq:Takesaki}
 \sigma_t^{\varphi_g}=g^{-1}\circ\sigma_t^\varphi\circ g.
 \end{equation}
 On the other hand, since $\varphi$ is $G$-strongly quasi invariant,
 \[
 \varphi_g(a)=\varphi(g(a))=\varphi(x_ga)
 \]
 for some positive element $x_g=x_g^*$.
 Therefore, for all $a\in \mathcal A$ and $g\in G$, since $x_g\in \mathrm{Centr}(\varphi)$,  we have \cite[Theorem 4.6]{PeTa}
 \begin{equation}\label{eq:phi_g_automorphism}
 \sigma_t^{\varphi_g}(a)=x_g^{it}\sigma_t^\varphi(a)x_g^{-it}.
 \end{equation}
 Combining \eqref{eq:Takesaki} and \eqref{eq:phi_g_automorphism}, we get \eqref{eq:G-strongly_quasi_invariance}.  \\
($\Longleftarrow$)  Suppose that \eqref{eq:G-strongly_quasi_invariance} holds:
\[
g^{-1}\circ \sigma_t^\varphi\circ g=x_g^{it}\sigma_t^\varphi x_g^{-it}.
\]
Note that the l.h.s. is the modular group for $\varphi_g=\varphi\circ g$. The r.h.s. is the modular group of $\psi(\cdot):=\varphi(x_g\cdot)$. Since $\varphi_g$ satisfies the modular condition with respect to $\{\sigma_t^{\varphi_g}\}$ it satisfies the modular condition with respect to $\{\sigma_t^\psi\}$. By \cite[Theorem 2.2]{Su}, there exists a unique positive injective operator $h_g$ affiliated with the center of $\mathcal A$ such that
\begin{equation}\label{eq:states}
\varphi_g(a)=\psi(h_ga)=\varphi(x_gh_ga), \quad \text{for all }a\in \mathcal A.
\end{equation}
Since $\mathcal A$ is a factor, its center is  $\mathbb C\unit$, and hence $h_g=\lambda_g\unit$ for some $\lambda_g>0$. By putting $y_g=\lambda_g x_g$, \eqref{eq:states} says that for all $g\in G$,
\[
\varphi_g(a)=\varphi(g(a))=\varphi(y_ga), \quad a\in \mathcal A,
\]
which proves that $\varphi$ is $G$-strongly quasi invariant.
\EndProof \\
We say that $(\sigma_t^\varphi)_{t\in \mathbb R}$ commutes with $G$ if $\sigma_t^\varphi$ commutes with $g$ for all $t\in \mathbb R$ and $g\in G$. The following theorem is an extension of \cite[Theorem 5.5]{PeTa} to the case of strong quasi invariance.
 \begin{thm}\label{thm:G-commuting_automorphism}
 The modular automorphism group $(\sigma_t^\varphi)_{t\in \mathbb R}$ commutes with $G$  if and only if $\varphi$ is $G$-strongly quasi invariant in the generalized sense and for each $g\in G$, the Radon-Nikodym derivative $x_g$ is affiliated with the center $\mathcal Z=\mathcal A\cap \mathcal A'$.
 \end{thm}
 \Proof
 ($\Longleftarrow$) Suppose that $\varphi$ is $G$-strongly quasi invariant in the generalized sense and the Radon-Nikodym derivative $x_g$ is affiliated with the center $\mathcal Z=\mathcal A\cap \mathcal A'$ for each $g\in G$.
Then for all $g\in G$, we have
\[
\varphi(g(a))=\varphi(x_ga), \quad \forall\,a\in \mathcal A.
\]
The modular group of $\varphi_g=\varphi\circ g$  is $g^{-1}\circ\sigma_t^\varphi\circ g$ and by the  uniqueness it  coincides with the modular group of the state $\varphi(x_g\cdot)$, which is $x_g^{it}\sigma_t^\varphi(\cdot)x_g^{-it}=\sigma_t^\varphi(\cdot)$ because $x_g^{\pm it}\in \mathcal Z$. We have seen that $\sigma_t^{\varphi}$ commutes with $G$.\\
 ($\Longrightarrow$) Suppose that $\{\sigma_t^\varphi\}$ commute with $G$. Then for each $g\in G$, we have
 \[
 \sigma_t^{\varphi_g}=g^{-1}\circ\sigma_t^\varphi\circ g=\sigma_t^\varphi.
 \]
Now, since $\sigma_t^{\varphi_g}=\sigma_t^\varphi$ is the automorphism group of $\varphi_g$, we conclude that $\varphi_g$ is $\sigma_t^\varphi$-invariant. Then, by Theorem 2.2 of \cite{Su}, there exists a unique positive injective $x_g$ affiliated with the center $\mathcal Z$ such that
 \[
 \varphi_g(a)=\varphi(g(a))=\varphi(x_g a),\quad a\in \mathcal A.
 \]
 This shows that $\varphi$ is $G$-strongly quasi invariant in the generalized sense with Radon-Nikodym derivatives $x_g$ being affiliated with the center $\mathcal Z$.
 \EndProof \\
If $\varphi$ is $G$-invariant, then obviously $\varphi$ is $G$-strongly quasi invariant in the generalized sense with Radon-Nikodym derivatives $x_g=\unit\in \mathcal Z$. Therefore we get
 \begin{cor}\label{cor:invariant_case}
 If $\varphi$ is $G$-invariant, then the modular automorphism group $(\sigma_t^\varphi)_{t\in \mathbb R}$ commutes with $G$.
 \end{cor}
As we will see in the example later, there are $G$-strongly quasi invariant states for which $\mathcal C$ is not a subalgebra of the center $\mathcal Z$. It is worth to notice the relationship between the modular automorphism group and the $G$-actions. Notice that $\mathrm{Centr}(\varphi)$ is equal to the invariant subalgebra for the modular automorphism group $(\sigma_t^\varphi)$ \cite{Ta}.
\begin{prop}\label{prop:G-strongly_quasi_G-action} Suppose that a faithful normal state $\varphi$ on a von Neumann algebra $\mathcal A$ is $G$-strongly quasi invariant. Let $(\sigma_t^\varphi)$ be the modular automorphism group of $\varphi$. Then, it holds that
\[
\varphi(g(\sigma_t^\varphi(a)))=\varphi(\sigma_t^\varphi(g(a))),\quad a\in\mathcal{A}.
\]
In particular, the state $\varphi_g=\varphi\circ g$ is $\sigma^\varphi$-invariant for all $g\in G$.
\end{prop}
 \Proof  The proof can be obtained by applying \cite[Theorem 2.4]{Su} (equivalence of (1) and (4) in that theorem), but for the self-containedness we provide a proof. We have
  \begin{eqnarray*}
  \varphi(g(\sigma_t^\varphi(a)))&=&\varphi(x_g\sigma_t^\varphi(a))\\
                                              &=&\varphi(\sigma_t^\varphi(x_ga))\quad (\text{invariance of }x_g\text{ under }\sigma^\varphi)\\
 &=&\varphi(x_ga)\quad (\sigma^\varphi\text{-invariance of }\varphi)\\
 &=&\varphi(g(a))\\
 &=&\varphi(\sigma_t^\varphi(g(a)))\quad (\sigma^\varphi\text{-invariance of }\varphi),
  \end{eqnarray*}
  showing the statement.
\EndProof \\
 Recall that  $\mathcal{F}(G)$ is the $G$-invariant subalgebra.
 \begin{lem}\label{lem:x_g_and_G-invariance}
 Suppose that $\varphi$ is a faithful $G$-strongly quasi invariant state on a von Neumann algebra $\mathcal A$. If any element $x_g$ belongs to $\mathcal F(G)$, then $x_g=\unit$.
 \end{lem}
 \Proof
Suppose that $x_{g}\in \mathcal{F}(G)$. Since $\mathcal{F}(G)$  ia a subalgebra of $\mathcal{A}$, $x_g^n\in\mathcal{F}(G)$ for all $n\in\mathbb{N}$. Thus, for all $n\ge 1$,
  \[
 \varphi(x_g^n)= \varphi(g(x_g^{n-1}))=\varphi(g(x_g)^{n-1})=\varphi(x_g^{n-1})=\cdots=\varphi(x_g)=\varphi(g(\unit))=1.
 \]
 Therefore,
 \[
 \varphi((x_g-\unit)^2)= \varphi(x_g^2-2x_g+\unit)=0,
  \]
  and we conclude $ x_g=\unit$ by the faithfulness of $ \varphi$.
   \EndProof
 \begin{thm}\label{thm:G-invariance}
 Let $\varphi$ be a faithful normal state on a von Neumann algebra $\mathcal A$ and let $G$ be a group of $*$-automorphisms of $\mathcal A$. Then the following conditions are equivalent:
 \begin{enumerate}
 \item[(i)] $\varphi$ is $G$-invariant\\[-3ex]
 \item[(ii)] $\varphi$ is $G$-strongly quasi invariant and $\mathcal C\sbs \mathcal F(G)$.\\[-3ex]
 \item[(iii)] $G$ commutes with $\sigma_t^\varphi$, $t\in \mathbb R$, (hence it is  $G$-strongly quasi invariant in the generalized sense) and for all $g\in G$, $ x_g\in \mathcal F(G)$.
 \end{enumerate}
 \end{thm}
 \Proof ((i)$\Longrightarrow$(ii)) Suppose (i) holds. Then, obviously for all $g\in G$ and $a\in \mathcal A$,
 \[
 \varphi(g(a))=\varphi(a).
 \]
 That is, $\varphi$ is $G$-strongly quasi invariant with $x_g\equiv \unit$. It follows that $\mathcal C\sbs \mathcal F(G)$. \\
 ((ii)$\Longrightarrow$(iii)) Suppose (ii). Since any element $x_g$ belongs to $\mathcal F(G)$, by Lemma \ref{lem:x_g_and_G-invariance} $x_g=\unit$. In particular $x_g$ is affiliated with the center $\mathcal Z$. Then,  by Theorem \ref{thm:G-commuting_automorphism},  $G$ commutes with $(\sigma_t^\varphi)$. \\
 ((iii)$\Longrightarrow$(i)) Suppose (iii).  Since $\sigma_t^\varphi$ commutes with $G$ for all $t\in \mathbb R$, by Theorem \ref{thm:G-commuting_automorphism} $\varphi$ is $G$-strongly quasi invariant in the generalized sense, and together with the property that $ x_g\in \mathcal F(G)\sbs\mathcal A$, $\varphi$ is $G$-strongly quasi invariant with $\mathcal C\sbs \mathcal F(G)$. By Lemma \ref{lem:x_g_and_G-invariance} again, $x_g=\unit$ for all $g\in G$, and hence $\varphi(g(a))=\varphi(a)$, proving that $\varphi$ is $G$-invariant.
 \EndProof
 \\
 \begin{rem}\label{rem:sufficient_condition_commutativity} (i) We notice that Theorem 1 in \cite{HeTa} follows as a corollary to  Theorem \ref{thm:G-invariance}. In fact, the result of Theorem \ref{thm:G-invariance} is stronger than that of \cite[Theorem 1]{HeTa} in the sense that in \cite[Theorem 1]{HeTa} they imposed the property that $G$ fixes the center $\mathcal Z$ elementwise.   \\
(ii) If $\varphi$ is $G$-strongly quasi invariant,   $\mathcal C$ is a subalgebra of $\mathrm{Centr}(\varphi)$ \cite{AcDh}. Therefore, the condition $\mathrm{Centr}(\varphi)\sbs\mathcal F(G)$ implies that $\mathcal C\sbs \mathcal F(G)$. Then, in that case by Theorem \ref{thm:G-invariance}, it follows that $\varphi$ is $G$-invariant and the group $G$ commutes with the modular automorphism group $(\sigma_t^\varphi)$.
\end{rem}
 
\subsection{Invariant subalgebras}
Let $\varphi$ be a $G$-strongly quasi invariant state on a von Neumann algebra $\mathcal A$.
In this subsection, we discuss the relationship between the invariant subalgebras of $G$ and of the modular automorphism group.  In the sequel we assume that $G$ is an amenable locally compact group. A group $G$ is said to be amenable if there exists an invariant mean over $G$, that is, a state $\eta$ on $\mathcal C_b(G)$, the space of bounded continuous functions on $G$, such that $\eta\{f(\hat g)\}=\eta\{f(h\hat g)\}=\eta\{f(\hat g h)\}$ for a fixed $h\in G$. Here $\hat g$ represents the dummy variable. Given a von Neumann algebra $\mathcal A$, as in \cite{DKKR}, in what follows we denote by $\mathcal A_{\mathcal C_b}^G$ the space of all weakly continuous, weakly$^*$ bounded functions from $G$ to $\mathcal A$ (particularly, for each $\Phi$ in the predual $\mathcal A_*$ of $\mathcal A$ and $X\in \mathcal A_{\mathcal C_b}^G$, $\Phi(X(\hat g))\in \mathcal C_b(G)$).  As in \cite[Lemma 1]{DKKR}, to a  mean $\eta$ over $G$ there exists a unique mapping $\tilde \eta$ from $\mathcal A_{\mathcal C_b}^G$ to  $(\mathcal A_*)^*=\mathcal A$ such that
\begin{equation}\label{eq:lift_mean}
\Phi(\tilde\eta(X))=\eta\{\Phi(X(\hat g))\}  \text{ for all }\Phi\in \mathcal A_*.
\end{equation}
Moreover, it
 satisfies that for all $a\in \mathcal A$,
\begin{equation}\label{eq:invariance_lift}
\tilde\eta(a)=a, \,\,\tilde\eta(a X)=a \tilde\eta(X), \,\,\tilde\eta(X a)=\tilde\eta(X) a,
\end{equation}
where $a$, $a X$, and $X a$ denote the elements of $\mathcal A_{\mathcal C_b}^G$ respectively defined by $a(g)=a$, $(a X)(g)=a X(g)$, and $(X a)(g)=X(g) a$. Given a mean $\eta$,  let us define a linear operator $\epsilon_\eta$ from  $\mathcal A$ to $\mathcal A$ by
\begin{equation}\label{eq:mean_operator}
\epsilon_\eta(x):=\tilde\eta(X_x),\quad x\in \mathcal A,
\end{equation}
where $X_x(g)=g(x)$, $g\in G$.
From the invariance of the mean, one easily checks that $\epsilon_\eta$ maps any element of $\mathcal A$ into the $G$-invariant subalgebra:
\begin{equation}\label{eq:group_invariant_projection}
g(\epsilon_\eta(x))=\epsilon_\eta(x),\quad x\in \mathcal A,\,\,g\in G.
\end{equation}
Also,  $\epsilon_\eta$ leaves $\mathcal F(G)$ elementwise invariant:
\begin{equation}\label{eq:elementwise_invariance}
\epsilon_\eta(x)=x,\quad x\in \mathcal F(G).
\end{equation}
\begin{prop}\label{prop:weaker_condition}
 Suppose that the function $g\mapsto \psi(g(x))$ is continuous for each   $\psi\in \mathcal A_*$ and let $\varphi$ be a faithful state on a von Neumann algebra $\mathcal A$ which is $G$-strongly quasi invariant. Then, for each mean $\eta$ on $G$, it holds that
\begin{equation}\label{eq:commutator_mean_modular}
\varphi\circ \epsilon_\eta=\varphi\circ\sigma_t^\varphi\circ \epsilon_\eta=\varphi\circ\epsilon_\eta\circ \sigma_t^\varphi
\end{equation}
\end{prop}
\Proof
We use Proposition \ref{prop:G-strongly_quasi_G-action} to see that for any $x\in \mathcal A$,
 \begin{eqnarray*}
 \varphi(\sigma_t^\varphi\circ \epsilon_\eta(x))&=&\varphi(\sigma_t^\varphi(\tilde\eta(X_x)))\\
 &=&\eta(\varphi(\sigma_t^\varphi\circ \hat g(x)))\\
 &=&\eta(\varphi(\hat g\circ\sigma_t^\varphi(x)))\\
  &=&\varphi(\tilde \eta(X_{\sigma_t^\varphi(x)}))\\
 &=&\varphi(\epsilon_\eta\circ \sigma_t^\varphi(x)).
  \end{eqnarray*}
This proves the statement.
\EndProof\\
Imposing a further condition, we get a stronger result.
  \begin{prop}\label{prop:mean_and_modular_group}
 Suppose that the function $g\mapsto \psi(g(x))$ is continuous for each  $\psi\in \mathcal A_*$ and let $\varphi$ be a faithful state on a von Neumann algebra $\mathcal A$ which is $G$-strongly quasi invariant.  If $G$ commutes with the automorphism group $(\sigma_t^\varphi)$, then for each mean $\eta$ on $G$, it holds that
\begin{equation}\label{eq:commutator_mean_modular}
\sigma_t^\varphi\circ \epsilon_\eta=\epsilon_\eta\circ \sigma_t^\varphi,\quad t\in \mathbb R.
\end{equation}
\end{prop}
\Proof
Let  $\psi\in \mathcal A_*$ be any  weakly continuous linear functional. For all $x\in \mathcal A$,
\begin{eqnarray*}
\psi(\sigma_t^\varphi(\epsilon_\eta(x)))&=&\psi(\sigma_t^\varphi(\tilde\eta(X_x)))\\
&=&\eta(\psi(\sigma_t^\varphi\circ\hat g(x)))\\
&=&\eta(\psi(\hat g\circ \sigma_t^\varphi(x)))  \quad( \text{commutativity of }\sigma_t^\varphi\text{ and } \hat g)\\
&=&\psi(\tilde \eta(X_{\sigma_t^\varphi(x)}))\\
&=&\psi(\epsilon_\eta(\sigma_t^\varphi(x))).
\end{eqnarray*}
Since $\psi$ is arbitrary, the statement follows.
\EndProof
\begin{cor}\label{cor:sufficient_condition}
Suppose that the function $g\mapsto \psi(g(x))$ is continuous for each  $\psi\in \mathcal A_*$ and let $\varphi$ be a faithful state on a von Neumann algebra $\mathcal A$ which is $G$-strongly quasi invariant. If $\mathcal C\sbs \mathcal Z$ or $\mathrm{Centr}(\varphi)\sbs\mathcal F(G)$, then \eqref{eq:commutator_mean_modular} holds.
\end{cor}
\Proof
 If  $\mathcal C\sbs \mathcal Z$, then by Theorem \ref{thm:G-commuting_automorphism}, $G$ commute with $(\sigma_t^\varphi)$. Hence by Proposition \ref{prop:mean_and_modular_group}, \eqref{eq:commutator_mean_modular} holds. For the latter case, the result follows from Remark \ref{rem:sufficient_condition_commutativity} (ii) and Proposition \ref{prop:mean_and_modular_group}.
\EndProof\\
Recalling that $\mathrm{Centr}(\varphi)$ is equal to the invariant subalgebra under the modular group $(\sigma_t^\varphi)$ we have
\begin{cor}\label{cor:inclusion} Under the conditions of Proposition \ref{prop:mean_and_modular_group}, it holds that for any mean $\eta$,
\begin{eqnarray*}
&&\epsilon_\eta(\mathrm{Centr}(\varphi))\sbs \mathrm{Centr}(\varphi),\\
&&\sigma_t^\varphi(\mathcal F(G))\sbs \mathcal F(G).
\end{eqnarray*}
\end{cor}
\Proof
If $x\in \mathrm{Centr}(\varphi)$, then $\sigma_t^\varphi(x)=x$ for all $t\in \mathbb R$. Thus,
\[
 \sigma_t^\varphi(\epsilon_\eta(x))=\epsilon_\eta(\sigma_t^\varphi(x))=\epsilon_\eta(x),\quad t\in \mathbb R,
\]
which implies that $\epsilon_\eta(x)\in \mathrm{Centr}(\varphi)$. The second relation follows immediately from the commutativity of $(\sigma_t^\varphi)$ and $G$.
\EndProof\\
Notice that $\mathbb R$ is an amenable group. Corresponding to $\epsilon_\eta$ defined in \eqref{eq:mean_operator}, we define a map $\tilde T$ from $\mathcal A$ to $\mathcal A$ for the automorphism group $(\sigma_t^\varphi)_{t\in \mathbb R}$.  Again, $\tilde T$ maps any element $x\in \mathcal A$ into the $\sigma_t^\varphi$-invariant subalgebra (cf. \eqref{eq:group_invariant_projection}): $\sigma_t^\varphi(\tilde T(x))=\tilde T(x)$ for any $t\in \mathbb R$ and $x\in \mathcal A$, implying that $\tilde T(x)\in \mathrm{Centr}(\varphi)$. One also notices that $\epsilon_\eta(x)\in \mathrm{co}[g(x)]^-$ and $\tilde T(x)\in \mathrm{co}[\sigma_t^\varphi(x)]^-$, the weakly closed convex hulls of $\{g(x):g\in G\}$ and $\{\sigma_t^\varphi(x):t\in \mathbb R\}$. Now we can state a main result in this subsection which corresponds to \cite[Theorem 3]{HeTa} with weaker conditions.
\begin{thm}\label{thm:ergodicity}
Let $\mathcal A$ be a von Neumann algebra and $G$ an amenable group. Let $\varphi$ be a faithful $G$-invariant state on $\mathcal A$.  Assume that the functions $g\mapsto \psi(g(x))$ and $t\mapsto \psi(\sigma_t^\varphi(x))$ are both continuous for each $\psi\in \mathcal A_*$. Suppose that
  $\mathrm{Centr}(\varphi)\sbs\mathcal F(G)$.
Then, $\epsilon_\eta(x)=\tilde T(x)$ for all $x\in \mathcal A$.
\end{thm}
\Proof
If $\mathrm{Centr}(\varphi)\sbs\mathcal F(G)$, then by Corollary \ref{cor:sufficient_condition},  the relation $\sigma_t^\varphi\circ \epsilon_\eta=\epsilon_\eta\circ \sigma_t^\varphi$ holds. Recall that $\epsilon_\eta$ and  $\tilde T$ are the projections on $\mathcal A$ with ranges $\mathcal F(G)$ and  $\mathrm{Centr}(\varphi)$, respectively, and  $\epsilon_\eta(x)\in \mathrm{co}[g(x)]^-$ and $\tilde T(x)\in \mathrm{co}[\sigma_t^\varphi(x)]^-$.  Proceeding now as in the proof of \cite[Theorem 3]{HeTa}, let $\sum_{i=1}^n\lambda_i\sigma_{t_i}^\varphi(x)\to \tilde T(x)$ strongly, $\sum_{i=1}^n\lambda_i=1$. We have by the continuity of $\epsilon_\eta$ that
\begin{equation}\label{eq:limit}
\epsilon_\eta\left(\sum_i^n\lambda_i\sigma_{t_i}^\varphi(x)\right)\to \epsilon_\eta\left(\tilde T(x)\right)=\tilde T(x),
\end{equation}
the equality holds because $\tilde T(x)\in \mathrm{Centr}(\varphi)\sbs\mathcal F(G)$ and $\epsilon_\eta$ leaves  $\mathcal F(G)$ elementwise invariant. On the other hand, one notes that $\epsilon_\eta\circ \sigma_t^\varphi(x)=\sigma_t^\varphi\circ\epsilon_\eta(x)=\epsilon_\eta(x)$. The second equality holds because $\sigma_t^\varphi$ leaves $\text{Centr}(\varphi)$, which is equal to the range of $\epsilon_\eta$, invariant elementwise. Therefore, the l.h.s. of \eqref{eq:limit} becomes
\[
\sum_i^n\lambda_i\epsilon_\eta(x)=\epsilon_\eta(x).
\]
We have shown the equality $\epsilon_\eta(x)=\tilde T(x)$ as desired.
\EndProof
\begin{rem}\label{rem:compare}
In \cite[Theorem 3]{HeTa}, the hypotheses of the theorem imply that  $\mathrm{Centr}(\varphi)=\mathcal F(G)=\mathcal Z$. However, in Theorem \ref{thm:ergodicity}, we can get rid of the condition of $\eta$-asymptotic abelianness.
\end{rem}
\subsection{Examples}
{\bf Example 1: Quasi invariant states}. \\[1ex]
We first give an example for a $G$-quasi invariant state following the idea given in \cite{AcDh}. Let $\mathcal A=\mathcal B(\mathcal H)$, the von Neumann algebra of all bounded linear operators on a separable Hilbert space $\mathcal H$. Let $\omega$ be a faithful normal state on $\mathcal A$ with a density operator $\rho\in \mathcal A$:
\begin{equation}\label{eq:normal_state}
\omega(a):=\tr{\rho a},\quad a\in \mathcal A.
\end{equation}
Let  $\sigma=(\sigma_t)_{t\in \mathbb R}$ be the group of modular automorphisms with respect to the state $\omega$:
\begin{equation}\label{eq:automorphism_group}
 \sigma_t(a):=\rho^{it}a\rho^{-it},\quad a\in \mathcal A,\,\,  t\in \mathbb R.
 \end{equation}
Let $K\in \mathcal A$ be any strictly positive operator with bounded inverse and normalized in the sense that
$\omega(K^{-1})=1$.
Let $G=\mathbb R$ and for each $t\in G$ define
 \begin{equation}\label{eq:RN_derivative1}
 x_t:=K\sigma_{-t}(K^{-1}).
 \end{equation}
 We define a state $\varphi=\varphi_K$ by
 \begin{equation}\label{eq:K-quasi_invariant_state1}
 \varphi(a):=\omega(K^{-1}a),\quad a\in \mathcal A.
 \end{equation}
 \begin{prop}\label{prop:K-quasi_invariant_state}
 The state $\varphi$ in \eqref{eq:K-quasi_invariant_state1} is $G$-quasi invariant with Radon-Nikodym derivatives given by \eqref{eq:RN_derivative1}.
 \end{prop}
 \Proof
 We see for all $t\in \mathbb R$ and $a\in \mathcal A$ that
 \begin{eqnarray*}
 \varphi(\sigma_t(a))&=&\omega(K^{-1}\sigma_t(a))\\
  &=&\omega(\sigma_{-t}(K^{-1})a)\quad (\sigma\text{-invariance of }\omega)\\
  &=&\varphi(K\sigma_{-t}(K^{-1})a)\\
  &=&\varphi(x_t a).
  \end{eqnarray*}
  It says that $\varphi$ is $G$-quasi invariant.
 \EndProof\\ [1ex]
 {\bf Example 2: Rotation group in $\mathcal M_2$}\\[1ex]
Let us find an example of $G$-strongly quasi invariant state in the simplest case.
Let $\mathcal A=\mathcal M_2(\mathbb C)$, the space of $2\times2$ matrices. Define
\begin{equation}\label{eq:normal_state}
\omega(a):=\frac{1}{2}\tr{a},\quad a\in \mathcal A.
\end{equation}
Let $G=\{g_\theta|\theta=0, \frac{\pi}{2}, \pi, \frac{3\pi}{2}\}$, where

 \begin{equation}\label{eq:normal_state}
g_\theta(a):=U_{-\theta} aU_{\theta},\quad a\in \mathcal A,\quad U_\theta=\left(\begin{matrix}\cos\theta&-\sin\theta\\\sin\theta&\cos\theta\end{matrix}\right).
\end{equation}
 Note that $G$ is a group of automorphisms on $\mathcal A$ and $\omega$ is a $G$-invariant state.

   Let $$\rho=\left(\begin{matrix}\lambda&0\\0&1-\lambda\end{matrix}\right),\quad \lambda=\frac{e^\beta}{1+e^\beta}\in(0,1).$$
 Define $K:=(2\rho)^{-1}$, which is a strictly positive operator with bounded inverse and it is normalized in the sense that
\begin{equation}\label{eq:normaization}
\omega({K^{-1}})=1.
\end{equation}
Let us define for each $\theta=0, \frac{\pi}{2}, \pi, \frac{3\pi}{2}$,
 \begin{equation}\label{eq:RN_derivative}
 x_{g_\theta}:=Kg_{-\theta}(K^{-1}).
 \end{equation}
Directly computing we get
 \begin{equation*}
x_{g_\theta}=(2\rho)^{-1}U_\theta(2\rho)U_{-\theta}=\rho^{-1}U_\theta \rho U_{-\theta}=\left(\begin{matrix}\frac{1+(2\lambda-1)\cos(2\theta)}{2\lambda}&\frac{2\lambda-1}{2\lambda}\sin(2\theta)\\\frac{2\lambda-1}{2(1-\lambda)}\sin(2\theta)&\frac{1+(1-2\lambda)\cos(2\theta)}{2(1-\lambda)}\end{matrix}\right).
 \end{equation*}
Concretely, we have
 \begin{equation}\label{eq:RN_derivatives}
 x_{g_0}=x_{g_{\pi}}=I,\quad x_{g_{{\pi}/{2}}}=x_{g_{{3\pi}/{2}}}=\left(\begin{matrix}\frac{1-\lambda}{\lambda}&0\\0&\frac{\lambda}{1-\lambda}\end{matrix}\right)=\left(\begin{matrix}e^{-\beta}&0\\0&e^\beta\end{matrix}\right).
 \end{equation}
We define a state $\varphi=\varphi_K$ by
 \begin{equation}\label{eq:K-quasi_invariant_state}
 \varphi(a):=\omega(K^{-1}a)=\tr{\rho a},\quad a\in \mathcal A.
 \end{equation}
 \begin{prop}\label{prop:K-G-strongly quasi_invariant_state}
 The state $\varphi$ in \eqref{eq:K-quasi_invariant_state} is a $G$-strongly quasi invariant state with Radon-Nikodym derivatives given by \eqref{eq:RN_derivative}.
 \end{prop}
 \Proof
 The proof of quasi invariance is the same as that of Theorem \ref{prop:K-quasi_invariant_state}. The positivity of the Radon-Nikodym derivatives are obvious from \eqref{eq:RN_derivatives}.
 \EndProof \\
Based on the above example, some remarks follow in relevance with the general theory developed in Section \ref{sec:general_theories}.
 \begin{rem}\label{rem:miscellany} Let $\varphi$ be the $G$-strongly quasi invariant state on $\mathcal A=\mathcal M_2(\mathbb C)$ considered in the above example.\\
 1. We compute that for $\theta=\pi/2$, $x_{g_\theta}=\left(\begin{matrix}\frac{1-\lambda}{\lambda}&0\\0&\frac{\lambda}{1-\lambda}\end{matrix}\right)$ and $g_{-\theta}(x_{g_\theta})=\left(\begin{matrix}\frac{\lambda}{1-\lambda}&0\\0&\frac{1-\lambda}{\lambda}\end{matrix}\right)\neq x_{g_\theta}$, showing that $\mathcal C\nsubseteq\mathcal F(G)$. By Theorem \ref{thm:G-invariance} $\varphi$ can't be $G$-invariant, as is actually the case.\\
	 2. One sees from \eqref{eq:RN_derivatives} that for $\theta=\pi/2$ and $3\pi/2$, $x_{g_\theta}$  do not belong to the center $\mathcal{Z}=\mathcal{A}\cap\mathcal{A}'$, which is $\mathbb C\unit$ in this model. By Theorem \ref{thm:G-commuting_automorphism}, $(\sigma_t^\varphi)_{t\in\mathbb{R}}$ do not commute with $G$, as one can directly check it. In fact, since $\sigma_t^{\varphi}(x)=\rho^{it}x\rho^{-it}$, for $\theta=\pi/2$ and $x= \left(\begin{matrix}a&b\\c&d\end{matrix}\right)\in \mathcal A$, one sees that
\[
\sigma_t^\varphi(g_\theta(x))=\left(\begin{matrix}d&-\left(\frac{\lambda}{1-\lambda}\right)^{it}c\\-\left(\frac{1-\lambda}{\lambda}\right)^{it}b&a\end{matrix}\right),\quad
g_\theta(\sigma_t^\varphi(x))=\left(\begin{matrix}d&-\left(\frac{1-\lambda}{\lambda}\right)^{it}c\\-\left(\frac{\lambda}{1-\lambda}\right)^{it}b&a\end{matrix}\right),
\]
showing that $\sigma_t^\varphi(g_\theta(x))\neq g_\theta(\sigma_t^\varphi(x))$. On the other hand, one checks that $\varphi(\sigma_t^\varphi(g_\theta(x)))=\varphi(g_\theta(\sigma_t^\varphi(x)))$ as Proposition \ref{prop:G-strongly_quasi_G-action} says.
 \end{rem}
 The Haar measure of the group in this example is the (normalized) uniform distribution. Therefore, we compute
\begin{equation}\label{eq:mean_ex1}
\kappa=\int_G x_gdg=\frac12\left(\begin{matrix}\frac1\lambda&0\\0&\frac1{1-\lambda}\end{matrix}\right).
\end{equation}
Hence, the $G$-invariant state $\varphi_G$ becomes
\begin{equation}\label{eq:invariant_state_ex1}
\varphi_G(a)=\varphi(\kappa a)=\tr{\rho \kappa a}=\frac12\tr{a}.
\end{equation}
The Umegaki conditional expectation is computed as
\begin{equation}\label{eq:umegaki_ex1}
E_G(a)=\frac12\left(\begin{matrix}a_{11}+a_{22}&a_{12}-a_{21}\\-(a_{12}-a_{21})&a_{11}+a_{22}\end{matrix}\right),\quad a=\left(\begin{matrix}a_{11}&a_{12}\\a_{21}&a_{22}\end{matrix}\right).
\end{equation}
Hence $E_G$ is the Umegaki conditional expectation onto the $G$-invariant subalgebra
\[
\mathcal F(G)=\left\{\left(\begin{matrix}a&b\\-b&a\end{matrix}\right):a,b\in \mathbb C\right\}.
\]
{\bf Example 3: Translation group in the cycles}\\[1ex]
 Let $C_N=\{0,1,\cdots,N-1\}$ be a cycle group of length $N$. Let $\mathcal A=\otimes_{i\in C_N}\mathcal M_2$, the $N$-tensor product of $\mathcal M_2$. Define an automorphism, cyclic translation, $\tau$ on $\mathcal A$ by
 \begin{equation}\label{eq:cyclic_translation}
 \tau(a)=a_1\otimes a_2\otimes\cdots a_{N-1}\otimes a_0\text { for }a=a_0\otimes a_1\otimes \cdots \otimes a_{N-1}\in \mathcal A.
 \end{equation}
 For each $n=0,1,\cdots,N-1$, let $g_n:=\tau^n$ and let $G=\{g_n|n=0,\cdots,N-1\}$. Notice that $g_n^{-1}=g_{N-n}$, $n=0,\cdots,N-1$.
Define a state $\omega$ on $\mathcal A$ by
\begin{equation}\label{eq:invariant_state_cycle}
\omega(a)=\frac1{2^N}\prod_{i=0}^{N-1}\tr{a_i},\quad  a=a_0\otimes a_1\otimes \cdots \otimes a_{N-1}\in \mathcal A.
\end{equation}
Obviously $\omega$ is $G$-invariant. Let $K:=\otimes_{i=0}^{N-1}K_i\in \mathcal A$ be a strictly positive operator with bounded inverse and assume that it is normalized in the sense that
\begin{equation}\label{eq:normaization_cycle}
\omega({K^{-1}})=1.
\end{equation}
For each $n\in \{0,\cdots,N-1\}$, let us define
 \begin{equation}\label{eq:RN_derivative_cycle}
 x_{g_n}:=Kg_{N-n}(K^{-1}).
 \end{equation}
 We assume that the operator $K$ satisfies $x_{g_n}=x_{g_n}^*$ for each $n=0,\cdots,N-1$. For instance, for the operator $K:=\otimes_{i=0}^{N-1}K_i$, if all the matrices $K_i$, $i=1,\cdots, N-1$, are diagonal and  strictly positive it satisfies the hypothesis. Now as in the previous examples define a state by
 \begin{equation}\label{eq:G-strongly_quasi_cycle}
  \varphi(a):=\omega(K^{-1}a),\quad a\in \mathcal A.
 \end{equation}
One checks that $\varphi$ is $G$-strongly quasi invariant with cocycles given in \eqref{eq:RN_derivative_cycle}.

The (normalized) Haar measure in this example is again a uniform distribution. Therefore we have
\[
\kappa=\frac1N \sum_{g\in G}x_g=\frac1N K\left(\sum_{n=0}^{N-1} g_n(K^{-1})\right).
\]
The $G$-invariant state $\varphi_G$ is then given by
\begin{eqnarray*}
\varphi_G(a)=\varphi(\kappa a)= \omega(K^{-1}\kappa a)&=&\frac1N \omega\left(\sum_{n=0}^{N-1}g_n(K^{-1})a\right) =\frac1N\omega\left(K^{-1}\sum_{n=0}^{N-1}g_n(a)\right)\\
&=&\frac1N\sum_{n=0}^{N-1}\varphi(g_n(a)),
\end{eqnarray*}
here in the second equality from the last we have used the fact that $\omega$ is $G$-invariant. Finally the Umegaki conditional expectation in this example is given by
\[
E_G(a)=\frac1N\sum_{n=0}^{N-1}g_n(a),\quad a\in \mathcal A.
\]

In particular, consider the case $K_i\equiv K_0$ for all $i=0,1,\cdots,N-1$. In this case $x_g=\unit$ for all $g\in G$ (see \eqref{eq:RN_derivative_cycle}) and $\varphi$ is $G$-invariant. Moreover, since the modular automorphism group of $\varphi$ is given by
\[
\sigma_t^\varphi(a)=K^{-it}aK^{it},
\]
one easily checks that $G$ and $(\sigma_t^\varphi)$ commute.\\[1ex]
{\bf Example 4: Spin flip dynamics}\\[1ex]
Here we provide with an example that satisfies the condition $\mathrm{Centr}(\varphi)\sbs \mathcal F(G)$, which was considered in Corollary \ref{cor:sufficient_condition} as well as in Theorem \ref{thm:ergodicity}. Let $\mathcal A=\mathcal M_2$, the space of $2\times 2$ matrices. Let $\sigma_z=\left(\begin{matrix}1&0\\0&-1\end{matrix}\right)$ be the $z$-component of the Pauli matrices. $\sigma_z$ is a unitary matrix: $\sigma_z^*=\sigma_z^{-1}=\sigma_z$, and we understand it also as an automorphism on $\mathcal A$ acting as (with a slight abuse of notation)
\begin{equation}\label{eq:flip_action}
\sigma_z(a)=\sigma_za\sigma_z^*=\left(\begin{matrix}a_{11}&-a_{12}\\-a_{21}&a_{22}\end{matrix}\right),\quad a=\left(\begin{matrix}a_{11}&a_{12}\\a_{21}&a_{22}\end{matrix}\right)\in \mathcal A.
\end{equation}
Then, $G=\{\unit,\sigma_z\}$ is a group of $*$-automorphisms on $\mathcal A$. Let $\rho=\left(\begin{matrix}\lambda&0\\0&\mu\end{matrix}\right)$ be a density matrix and define
a state $\varphi$ by
\[
\varphi(a):=\tr{\rho a},\quad a\in \mathcal A.
\]
One easily checks that $\varphi$ is $G$-invariant, and hence $G$-strongly quasi invariant. Also from the group action in \eqref{eq:flip_action}, one sees that $\mathcal F(G)$, the subalgebra of $G$-fixed elements, is the commutative subalgebra of $\mathcal A$ consisting of diagonal matrices. On the other hand, the modular automorphism of $\varphi$ is given by
\begin{equation}\label{eq:modular_group_spin_flip}
\sigma_t^\varphi(a)=\rho^{it}a\rho^{-it}= \left(\begin{matrix}a_{11}&(\lambda/\mu)^{it}a_{12}\\(\lambda/\mu)^{-it}a_{21}&a_{22}\end{matrix}\right),\quad  a=\left(\begin{matrix}a_{11}&a_{12}\\a_{21}&a_{22}\end{matrix}\right)\in \mathcal A.
\end{equation}
First as can be  easily checked from \eqref{eq:flip_action} and \eqref{eq:modular_group_spin_flip}, we notice that the group $G$ and modular automorphism group $(\sigma_t^\varphi)$ commute in this example. Moreover, if $\lambda\neq\mu$, the subalgebra of elements that are fixed by the modular group $(\sigma_t^\varphi)$, or $\mathrm{Centr}(\varphi)$ is also the subalgebra of diagonal matrices. We have shown
\begin{prop}\label{prop:spin_flip}
Let $\mathcal A=\mathcal M_2$ and let $\varphi$ and $G$ be a state and a group of $*$-automorphisms, respectively, defined above. Then, the group $G$ and the modular automorphism group $(\sigma_t^\varphi)$ commute. Moreover, if $\lambda\neq\mu$, then $\mathrm{Centr}(\varphi)=\mathcal F(G)$ and both are equal to the commutative subalgebra of diagonal matrices.
\end{prop}
\subsection{Classical spin systems}\label{subsec:spin_systems}
In this subsection we recall an example of Gibbs measures for classical spin systems. We refer to \cite{DKY} for the details. First we remark that for the definition of quasi-invariance we may consider the states on the $C^*$-algebras \cite{AcDh}.

We let $\Omega:=\{-1,1\}^{\mathbb Z^d}$ be the set of spin configurations on the integer lattice $\mathbb Z^d$. Let $\Phi=(\Phi_\Lambda)_{\Lambda\sbs\sbs\mathbb Z^d}$ be a translation invariant, finite range interaction, and $\mu$ a Gibbs measure for it \cite{Ge}. The most simple and well-known example is the Ising model:  for $\xi=(\xi_i)_{i\in \mathbb Z^d}\in \Omega$,
\[
\Phi_\Lambda(\xi)=\begin{cases}-J\xi_i\xi_j,&\Lambda=\{i,j\},\,\,|i-j|=1,\\ h\xi_i,&\Lambda=\{i\},\\
0,&\text{otherwise}\end{cases}.
\]
Here $J$ is the interaction strength and $h$ denotes the external magnetic field strength; $J>0$ for ferro magnetic model and $J<0$ for anti-ferro magnetic model. We can think of $\mu$ as a state on the $C^*$-algebra $\mathcal A:=C(\Omega)$, the space of continuous functions on the set $\Omega$.

For each $N\in \mathbb N$, let $\Lambda_N=[-N,N]^d\cap \mathbb Z^d$ denote the rectangular box with sides of length $2N+1$. Let $(G_N)_{N\in \mathbb N}$ be an increasing sequence of automorphisms of $\mathcal A$ such that for each $N\in \mathbb N$, $G_N$ depends only on the local configurations in $\Lambda_N$. We let $G=\cup_{N\in \mathbb N}G_N$. We consider $G$ of spin interchanges or spin flips defined as follows:
\begin{ex}\label{ex:spin_automorphism_group}
Notice that any continuous bijection $\tau:\Omega\to \Omega$ naturally induces an automorphism $\tau:\mathcal A\to \mathcal A$ by 
\[
\tau(f)(\omega)=f(\tau(\omega)),\quad f\in \mathcal A.
\]
\begin{enumerate} 
\item[(i)]  (Spin exchanges) For $i\neq j\in \mathbb Z^d$, $\tau_{ij}:\Omega\to \Omega$ is defined by
\[(\tau_{ij}(\omega))_k=\omega_k^{ij}:=\begin{cases}\omega_k,&k\neq i,j\\ \omega_j,&k=i,\\
\omega_i, &k=j.\end{cases}
\]
\end{enumerate}
The group $G_N$ is generated by $\{\tau_{ij}:i\neq j\in \Lambda_N\}$. In other words, $G_N$ consists of spin permutations in the box $\Lambda_N$.
\item[(ii)] (Spin flips) For each $i\in \mathbb Z^d$, $\tau_i:\Omega\to \Omega$ is defined by
\[(\tau_i(\omega))_j=\omega_j^i:=\begin{cases}\omega_j,&j\neq i,\\ -\omega_i,&j=i.\end{cases}  
\]
The group $G_N$ is generated by $\{\tau_i:i\in \Lambda_N\}$ so it is the group of partial spin flips in $\Lambda_N$. 
\end{ex}
The following was shown in \cite{DKY}.
\begin{thm}\label{thm:quasi_invariance_spin_classical}
Let $\Phi$ be a translation invariant and finite range interaction for the spin system and let $\mu$ be a Gibbs measure for $\Phi$. Let $G=\cup_{N\in \mathbb N}G_N$ be one of the locally compact groups introduced in Example \ref{ex:spin_automorphism_group}. Then $\mu$ is $G$-strongly quasi-invariant with cocycles $x_\tau$ given by 
\begin{equation}\label{eq:cocycles_spin}
x_\tau(\omega)=\exp[H(\omega)-H(\tau^{-1}(\omega))], \quad \tau\in G,
\end{equation}
here, the exponent is defined by 
\[
H(\omega)-H(\tau^{-1}(\omega))=\lim_{N\to \infty} \sum_{X\sbs \Lambda_N}\left(\Phi_X(\omega)-\Phi_X(\tau^{-1}(\omega))\right)\Big], 
\]
which is well-defined since $\tau$ gives only a local change. 
\end{thm} 
  
\section{$G$-strongly quasi invariant states and tracial states}\label{sec:tracial}
In this section we characterize the $G$-strongly quasi invariant states by tracial states with suitable densities. It is an analogous result to the one developed by St{\o}rmer \cite[Theorem 1]{St}. Given a semifinite von Neumann algebra $\mathcal A$, a group $G$ of $*$-automorphisms is said to acts ergodicly on the center $\mathcal Z$ of $\mathcal A$ if $\mathcal F(G)\cap \mathcal Z=\mathbb C \unit$, where $\mathcal F(G)$ is the fixed points of $G$ in $\mathcal A$.

Let $\mathcal A$ be a semifinite von Neumann algebra. Let $G$ be a locally compact amenable group of $*$-automorphisms acting ergodicly on the center $\mathcal Z$ of $\mathcal A$.  Let $\varphi$ be a faithful $G$-strongly quasi invariant state on $\mathcal A$. Let $\eta$ be a mean. \\[1ex]
{\bf Assumption}. We suppose the following conditions:
\begin{enumerate}
\item[-] the function $X:G\to \mathcal A$ defined by $X(g)=x_g$ is weakly continuous.
\item[-] $\tilde\eta(X)$ belongs to $\mathcal A$ with a bounded inverse $\tilde \eta(X)^{-1}$.
\end{enumerate}
\begin{lem}\label{lem:mean_density}
It holds that
\begin{enumerate}
\item[(i)] The functional $\varphi_G$ on $\mathcal A$ defined by $\varphi_G(a)=\varphi(\tilde\eta(X)a)$, $a\in \mathcal A$, is a $G$-invariant state.\\[-3ex]
\item[(ii)]  $\tilde\eta(X)$ is self-adjoint, positive definite, and belongs to $\mathrm{Centr}(\varphi)$.
\end{enumerate}
\end{lem}
\Proof
(i) Let us first show that $\varphi_G$ is a state on $\mathcal A$. For any $a\in \mathcal A$ we have
\begin{eqnarray}\label{eq:G_invariant_state_amenable}
\varphi_G(a)&=&\varphi(\tilde\eta(X)a)\nonumber\\
&=&\eta(\varphi(X(\hat g)a))\nonumber\\
&=&\eta(\varphi(x_{\hat g}a))\nonumber\\
&=&\eta(\varphi(\hat g(a))).
\end{eqnarray}
From this we see that $\varphi_G(\unit)=1$ and $\varphi_G(a)\ge 0$ for $a\ge 0$. Therefore, $\varphi_G$ is a state on $\mathcal A$. For $G$-invariance, we see from \eqref{eq:G_invariant_state_amenable} that for any $h\in G$,
\[
\varphi_G(h(a))=\eta(\varphi(\hat gh(a)))=\eta(\varphi(\hat g(a)))=\varphi_G(a).
\]
The second equality follows from the invariance of the mean.\\
(ii) For any $a\ge 0$,
\begin{equation}\label{eq:equality1}
\varphi_G(a)=\overline{\varphi_G(a)}=\overline{\varphi(\tilde\eta(X)a)}=\varphi(a\tilde\eta(X)^*).
\end{equation}
On the other hand, since $\varphi$ is $G$-strongly quasi invariant, from \eqref{eq:G_invariant_state_amenable} we get
\begin{equation}\label{eq:equality2}
\varphi_G(a)=\eta(\varphi(x_{\hat g}a))=\eta(\varphi(ax_{\hat g}))=\varphi(a\tilde\eta(X)).
\end{equation}
From \eqref{eq:equality1} and \eqref{eq:equality2} we have $\varphi(a\tilde\eta(X)^*)=\varphi(a\tilde\eta(X))$ for all $a\ge 0$. Now by using the spectral decomposition and the faithfulness of $\varphi$ we conclude that $\tilde\eta(X)$ is self-adjoint. To see the positive definiteness, let us decompose $\tilde\eta(X)=\tilde\eta(X)^+-\tilde\eta(X)^-$ into positive and negative parts. Then, since $\tilde\eta(X)^{\pm}$ are positive definite and orthogonal to each other, we have
\[
0\le \varphi_G(\tilde\eta(X)^-)=\varphi(\tilde\eta(X)\tilde\eta(X)^-)=-\varphi((\tilde\eta(X)^-)^2)\le 0.
\]
Since $\varphi$ is faithful we conclude that $\tilde\eta(X)^-=0$ and hence $\tilde\eta(X)=\tilde\eta(X)^+\ge 0$. Finally from \eqref{eq:G_invariant_state_amenable} we have for all $a\in \mathcal A$,
\[
\varphi(\tilde\eta(X)a)=\eta(\varphi(x_{\hat g}a))=\eta(\varphi(ax_{\hat g}))=\varphi(a\tilde\eta(X)).
\]
It shows that $\tilde\eta(X)$ belongs to $\mathrm{Centr}(\varphi)$.
\EndProof
\begin{thm}\label{thm:characterization_by_tracial_states_amenable_group}
Let $\mathcal A$ be a semifinite von Neumann algebra and let $G$ be an amenable group of  $*$-automorphisms of $\mathcal A$, and we assume that $G$ acts ergodicly on the center $\mathcal Z$ of $\mathcal A$. Suppose that $\varphi$ is a faithful $G$-strongly quasi invariant state on $\mathcal A$ with Radon-Nikodym derivatives $x_g$'s. Further, we suppose the conditions in the Assumption of this section. Then, there exists up to a scalar multiple a unique faithful normal $G$-invariant semifinite trace $\tau$ of $\mathcal A$, and there is a positive self-adjoint operator $c$ affiliated with $\rm{Centr}(\varphi)$ such that
$\varphi(a)=\tau(c a)$ for all $a\in \mathcal A$ and it holds that $g^{-1}(c)=cx_g$ for all $g\in G$.
\end{thm}
\Proof
In Lemma \ref{lem:mean_density}  we have seen that the state $\varphi_G(\cdot)=\varphi(\tilde\eta(X)\cdot)$ is $G$-invariant. Applying \cite[Theorem 1]{St} for $\varphi_G$, there exists a faithful normal $G$-invariant semifinite trace $\tau$ and a positive operator $b$ affiliated with $\mathcal F(G)$ such that $\varphi_G(a)=\tau(ba)$, and hence $\varphi(a)=\tau(b\tilde\eta(X)^{-1}a)$ for all $a\in \mathcal A$.
If $b$ and $\tilde\eta(X)^{-1}$ commute, then $c:=b\tilde\eta(X)^{-1}$ is positive and we are done since $c=b\tilde\eta(X)^{-1}$ is also affiliated with $\mathcal F(G)$. So, it remains to show that $b$ and $\tilde\eta(X)^{-1}$ commute. Notice that if $a\ge 0$, then $\varphi(a)=\tau(b\tilde\eta(X)^{-1}a)\ge 0$. This implies that
\[
\tau(b\tilde\eta(X)^{-1}a)=\overline{\tau(b\tilde\eta(X)^{-1}a)}=\tau((b\tilde\eta(X)^{-1}a)^*)=\tau(a\tilde\eta(X)^{-1}b)=\tau(\tilde\eta(X)^{-1}ba).
\]
This shows that $\tau([b,\tilde\eta(X)^{-1}]a)=0$ for all $a\ge 0$. Applying polar decomposition for any $a\in \mathcal A$, we conclude that $\tau([b,\tilde\eta(X)^{-1}]a)=0$ for all $a \in \mathcal A$, showing that $[b,\tilde\eta(X)^{-1}]=0$.

To show the last relation, we see that for any $a\in \mathcal A$ and $g\in G$
\begin{eqnarray*}
\tau(g^{-1}(c)a)&=&\tau(cg(a))\quad(G\text{-invariance of }\tau)\\
&=&\varphi(g(a))\\
&=&\varphi(x_ga)\\
&=&\tau(cx_ga).
\end{eqnarray*}
It proves the statement.
\EndProof
\begin{rem} \label{rem:compact}
Suppose that $G$ is a compact group acting ergodicly on a semifinite von Neumann algebra $\mathcal A$ and let $\varphi$ be a $G$-strongly quasi invariant on $\mathcal A$. We suppose that the map  $g\mapsto x_g$ is continuous. Considering the mean $\eta$ as the integral w.r.t. the Haar measure, we see that $\tilde \eta(X)$ is nothing but $\kappa=\int_Gx_gdg$ in \eqref{eq:rnd_mean_and_Umegaki}. Since $\kappa\in \mathcal A$ with a bounded inverse $\kappa^{-1}$, the conditions of the Assumption are fulfilled.
\end{rem}

\section{Modular theory for strongly quasi invariant states}\label{sec:modular}
In this section we discuss an application of strongly quasi invariant states. Let $\mathfrak M$ be a von Neumann algebra on a separable Hilbert space $\mathfrak h$ and $\Omega$ a cyclic and separating vector for $\mathfrak M$. Considering a natural positive cone $\mathcal P$ associated with the pair $(\mathfrak M, \Omega)$, there is a well known theory for the relation of automorphism group of $\mathfrak M$ and standard forms defined by the elements of $\mathcal P$ (see Section 2.5.4 of \cite{BR}). In this section we show that the parallel theory can be developed when we consider the strongly quasi invariant states.

Throughout this section we fix a $G$-strongly quasi invariant state $\varphi$ on a von Neumann algebra $\mathcal A$. Let $(\mathcal H_\varphi,\pi_\varphi,\Omega_\varphi)$ be the cyclic representation associated with the state $\varphi$ \cite{BR}. For each $g\in G$, let $\varphi_g$ be the state on $\mathcal A$ given by
 \[
 \varphi_g(a)=\varphi(g(a)),\quad a\in \mathcal A.
 \]
Since $\varphi$ is $G$-strongly quasi invariant, we have
\begin{eqnarray}\label{eq:g-GNS}
 \varphi_g(a)=\varphi(x_ga)&=&\varphi(\sqrt{x_g}a\sqrt{x_g})\nonumber\\
 &=&\langle \pi_\varphi(\sqrt{x_g})\Omega_\varphi,\pi_\varphi(a)\pi_\varphi(\sqrt{x_g})\Omega_\varphi\rangle.
 \end{eqnarray}
Defining $\mathcal H:=\mathcal H_\varphi$,  $\pi:=\pi_\varphi$, and $\Omega_g=\pi_\varphi(\sqrt{x_g})\Omega_\varphi$, we conclude from \eqref{eq:g-GNS} that $(\mathcal H,\pi,\Omega_g)$ is a cyclic representation associated with $\varphi_g$.

Now let $J_\varphi$ and $\Delta_\varphi$ be the modular conjugation and modular operator, respectively, corresponding to the cyclic and separating vector $\Omega_\varphi$. Similarly, given a group element $g\in G$ let $J_g$ and $\Delta_g$ be the modular conjugation and modular operator corresponding to the cyclic and separating vector $\Omega_g$. We would like to find the relations of the two systems. Let $\mathcal M:=\pi(\mathcal A)''$ and recall that a natural positive cone $\mathcal P=\mathcal P_\varphi$ associated with the pair $(\mathcal M,\Omega_\varphi)$ is the closure of the set (see \cite[Definition 2.5.25]{BR})
\[
\{\pi(a)j(\pi(a))\Omega_\varphi:a\in \mathcal A\},
\]
where $j(\pi(a)):=J_\varphi\pi(a)J_\varphi$. First we observe:
\begin{prop}\label{prop:unitary_relation}
For each $g\in G$, we have the following properties:
\begin{enumerate}
\item[(i)]  $\Omega_g$ belongs to the natural positive cone $\mathcal P$ associated with $\Omega_\varphi$.\\[-3ex]
 \item[(ii)] $J_\varphi=J_g$ and it follows that $J_\varphi\Omega_g=\Omega_g$.
 \end{enumerate}
 \end{prop}
 \Proof (i) Since $\pi(x_g)$ commutes with $\Delta_\varphi$,
 \begin{equation}\label{eq:positive vector}
 \Omega_g=\pi (\sqrt{x_g})\Omega_\varphi=\pi (\sqrt{x_g})\Delta_\varphi^{1/4}\Omega_\varphi=\Delta_\varphi^{1/4}\pi (\sqrt{x_g})\Omega_\varphi.
 \end{equation}
 By \cite[Proposition 2.5.26]{BR}, $\Omega_g$ belongs to $\mathcal P$. \\
   (ii) Since the cyclic and separating vector $\Omega_g$ belongs to the natural positive cone associate with $\Omega_\varphi$, by \cite[Proposition 2.5.30]{BR} it follows that $J_\varphi=J_g$ and so $J_\varphi\Omega_g=J_g\Omega_g=\Omega_g$.
    \EndProof  \\
 From now on we write $J:=J_\varphi=J_g$.
\begin{prop}\label{prop:unitary_relation}
For each $g\in G$, the operaor
\begin{equation}\label{eq:unitary_map}
 U_g: \pi(a)\Omega_g\mapsto \pi(g(a))\Omega_\varphi,\quad a\in \mathcal A,
 \end{equation}
extends to a unitary operator and the relation holds:
 \begin{equation}\label{eq:unitarity_for_quasi_invariant_state}
u_g(\pi(a)):= U_g\pi(a)U_g^{-1}=\pi(g(a)), \quad a\in \mathcal A.
 \end{equation}
 \end{prop}
 \begin{rem}\label{rem:unitarity_for_quasi_invariant_state} (i) Given an automorphism $g$, it is a standard way to define a unitary operator $U_g$ as \eqref{eq:unitary_map} which satisfies the relation \eqref{eq:unitarity_for_quasi_invariant_state} \cite[Section 4]{Su}, \cite[Corollary 2.5.32]{BR}. We emphasize, however, that from the property of $G$-strongly quasi invariance, $\{U_g:g\in G\}$ becomes a group as noted below and thereby $\{u_g:g\in G\}$ is a representation of $G$ on $\mathcal M$.

(ii) The relation \eqref{eq:unitarity_for_quasi_invariant_state} is developed in the case of invariant states (\cite[Corollary 2.3.17]{BR}) and so the above result is an extension to the $G$-strongly quasi invariant case. The unitary operator $U_g$ can also be defined as
 \begin{equation}\label{eq:unitary_other_representation}
 U_g(\pi(a)\Omega_\varphi)=\pi(g(a)x_{g^{-1}}^{1/2})\Omega_\varphi.
 \end{equation}
 In fact, from $\Omega_g=\pi(x_g^{1/2})\Omega_\varphi$,
 \begin{eqnarray*}
 U_g(\pi(a)\Omega_\varphi)&=&U_g(\pi(ax_g^{-1/2})\Omega_g)\\
 &=&\pi(g(a)g(x_g^{-1/2}))\Omega_\varphi\\
 &=&\pi(g(a)x_{g^{-1}}^{1/2})\Omega_\varphi.\quad(g(x_g^{-1})=x_{g^{-1}})
 \end{eqnarray*}
 The definition \eqref{eq:unitary_other_representation} was introduced in \cite{AcDh}. Moreover, in \cite[Theorem 6]{AcDh}, it was shown that $\{U_g\}_g$ is a unitary representation of $G$, particularly meaning that
 \[
 U_gU_h=U_{gh}\text{ and }U_g^*=U_{g^{-1}},\quad g,h\in G.
 \]
 \end{rem}
 \Proof[ of Proposition \ref{prop:unitary_relation}]
Notice that $\{\pi(a)\Omega_\varphi:a\in \mathcal A\}$ and $\{\pi(a)\Omega_g:a\in \mathcal A\}$ are both dense in $\mathcal H$.  By using the cyclic representation associated with $\varphi_g$ introduced above, for any $a,b\in \mathcal A$,
 \begin{eqnarray*}
 \langle \pi(a)\Omega_g,\pi(b)\Omega_g\rangle&=& \langle  \Omega_g,\pi(a^*b)\Omega_g\rangle\\
 &=&\varphi_g(a^*b)\\
 &=&\varphi(g(a^*b))=\varphi(g(a)^*g(b))\\
 &=&\langle\Omega_\varphi,\pi(g(a)^*g(b))\Omega_\varphi\rangle\\
 &=&\langle \pi(g(a))\Omega_\varphi,\pi(g(b))\Omega_\varphi\rangle.
 \end{eqnarray*}
 Therefore, the map $U_g$ extends to a unitary operator on $\mathcal H$. As $\{\pi(a)\Omega_g:a\in \mathcal A\}$ is dense in $\mathcal H$, the relation \eqref{eq:unitarity_for_quasi_invariant_state} is equivalent to
 \begin{equation}\label{eq:equivalent_relation}
 U_g\pi(a)\pi(b)\Omega_g=\pi(g(a))U_g\pi(b)\Omega_g,\quad a,b\in \mathcal A.
 \end{equation}
In fact,
\[
\text{ l.h.s. of \eqref{eq:equivalent_relation}}=U_g\pi(ab)\Omega_g= \pi(g(ab))\Omega_\varphi,
\]
 and
 \[
 \text{r.h.s. of \eqref{eq:equivalent_relation}}=\pi(g(a))\pi(g(b))\Omega_\varphi=\pi(g(ab))\Omega_\varphi,
 \]
 proving the equality of \eqref{eq:equivalent_relation}.
 \EndProof\\
Recall \cite{BR} that
\[
S_\varphi=J\Delta_\varphi^{1/2} \text{ and }S_g=J\Delta_g^{1/2},
\]
where $S_\varphi$ and $S_g$ are the (closures of ) operators defined on $\mathcal M \Omega_\varphi=\mathcal M\Omega_g$ by
\[
S_\varphi(\pi(a)\Omega_\varphi)=\pi(a^*)\Omega_\varphi\text{ and } S_g(\pi(a)\Omega_g)=\pi(a^*)\Omega_g.
\]
Similarly, recalling the relation $J\mathcal{M}J=\mathcal{M}'$, the commutant of $\mathcal{M}$, the operators $F_\varphi$ and $F_g$ are defined on $\mathcal M'\Omega_\varphi$ and $\mathcal M'\Omega_g$, respectively, by
\[
F_\varphi(J\pi(a)J\Omega_\varphi)=J\pi(a^*)J\Omega_\varphi \text{ and }F_g(J\pi(a)J\Omega_g)=J\pi(a^*)J\Omega_g.
\]
\begin{prop}\label{prop:anti-linear_operators}
The following relation holds:
\begin{equation}\label{eq:exchange_of_operators}
S_\varphi U_g=U_gS_g \text{ on }\pi(\mathcal A)\Omega_g.
\end{equation}
Furthermore, the following relations hold:
\begin{equation}\label{eq:further commuting}
 JU_g=U_gJ\text{ and }\Delta_\varphi U_g=U_g\Delta_g \text{ on }\pi(\mathcal A)\Omega_g.
 \end{equation}
\end{prop}
\Proof
For any $\pi(a)\Omega_g\in \pi(\mathcal A)\Omega_g$, it holds
\begin{eqnarray*}
S_\varphi U_g(\pi(a)\Omega_g)&=&S_\varphi(\pi(g(a))\Omega_\varphi)\\
&=&\pi(g(a^*))\Omega_\varphi\\
 &=&U_g(\pi(a^*)\Omega_g)\\
 &=&U_gS_g(\pi(a)\Omega_g).
\end{eqnarray*}
This proves \eqref{eq:exchange_of_operators}. To show \eqref{eq:further commuting}, we see  from \eqref{eq:exchange_of_operators} that on $\pi(\mathcal A)\Omega_\varphi$
\begin{eqnarray*}
S_\varphi&=&U_gS_gU_g^*\\
&=&U_gJ\Delta_g^{1/2}U_g^*\\
&=& U_gJU_g^*U_g\Delta_g^{1/2}U_g^*.
\end{eqnarray*}
Now the polar decomposition of $S_\varphi=J\Delta_\varphi^{1/2}$ and applying the uniqueness of polar decomposition the result follows.
\EndProof
 \begin{rem} (i) The unitary operator $U_g$ is represented on another dense subspace $J\pi(\mathcal A)J\Omega_g$ as follows:
 \begin{eqnarray*}
 U_gJ\pi(a)J\Omega_g&=& U_gJ\pi(a)\Omega_g\\
&=& JU_g\pi(a)\Omega_g\\
 &=&J\pi(g(a))\Omega_\varphi=J\pi(g(a))J\Omega_\varphi.
 \end{eqnarray*}
 (ii) By using (i) above and the same method used in the proof of \eqref{eq:exchange_of_operators}, we can show that
 \[
 F_\varphi U_g=U_g F_g\text{ on }J\pi(\mathcal A)J\Omega_g.
 \]
 (iii) The result of Proposition \ref{prop:anti-linear_operators} corresponds to \cite[Corollary 2.5.32]{BR} and \cite[Lemma 2]{St} for the case of strongly quasi invariant states.
 \end{rem}
 \begin{thm}\label{thm:modular operator relation}
 The following two equivalent relations hold:
 \begin{enumerate}
 \item[(i)] $S_g=\pi(\sqrt{x_g^{-1}})J\pi(\sqrt{x_g})JS_\varphi$
 \item[(ii)] $\Delta_g^{1/2}=\pi(\sqrt{x_g})J\pi(\sqrt{x_g^{-1}})J\Delta_\varphi^{1/2}$
 \end{enumerate}
 \end{thm}
 \Proof The equivalence of (i) and (ii) is clear. To show (i), since $\Omega_g=J\Omega_g=J\pi(\sqrt{x_g})\Omega_\varphi=J\pi(\sqrt{x_g})J\Omega_\varphi$, for any $a\in\mathcal{A}$,
 \begin{eqnarray*}
 S_g\pi(a)\Omega_g&=&\pi(a^*)\Omega_g\\
 &=&\pi(a^*)J\pi(\sqrt{x_g})J\Omega_\varphi\\
 &=&J\pi(\sqrt{x_g})J\pi(a^*)\Omega_\varphi\\
 &=&J\pi(\sqrt{x_g})J\pi(\sqrt{x_g^{-1}})\pi(\sqrt{x_g}a^*)\Omega_\varphi\\
 &=&J\pi(\sqrt{x_g})J\pi(\sqrt{x_g^{-1}})S_\varphi\pi(a\sqrt{x_g})\Omega_\varphi\\
 &=&\pi(\sqrt{x_g^{-1}})J\pi(\sqrt{x_g})JS_\varphi\pi(a)\Omega_g,
  \end{eqnarray*}
  which implies the equality.
 \EndProof

\section{Projection onto the invariant subspace}\label{sec:abelian}
In this section we investigate the orthogonal projection onto the invariant subspace for the group $\{U_g:g\in G\}$ defined in the previous section. We start by discussing the abelian subalgebras. The related results can be found in \cite{DKKR, DKS}.

Let us suppose that $G$ is a compact group. In this case define
\begin{equation}\label{eq:projection}
P_G:=\int_G U_gdg.
\end{equation}
Since the functions $g\mapsto g(a)$, $a\in \mathcal A$, and $g\mapsto x_g$ are assumed to be strongly continuous in Section \ref{sec:general_theories}, the above integral is well defined. It can be shown that \cite{AcDh} $P_G$ is a projection onto the closure of
\begin{equation}\label{eq:projection_range}
\{\psi\in \mathcal H_\varphi:U_g(\psi)=\psi,\,\,\,\forall\,g \in G\}.
\end{equation}
In fact, by \eqref{eq:projection} we have for any $g\in G$,
\[
P_GU_g=U_gP_G=P_G.
\]
Let $\mathcal R$ be the von Neumann algebra generated by $\pi(\mathcal A)$ and $\{U_g:g\in G\}$.
By integrating both sides of \eqref{eq:unitarity_for_quasi_invariant_state} we get
\[
\int_Gu_g(\pi(a))dg=\int_G\pi(g(a))dg=\pi(E_G(a))=:\tilde E_G(\pi(a)).
\]
\begin{rem}\label{rem:Umegaki} For any $a\in \mathcal A$, $\tilde E_G(\pi(a))$ is invariant under $u_g$ for all $g\in G$. In fact, for any $a\in \mathcal A$
\begin{eqnarray*}
u_g \tilde E_G(\pi(a))&=& \int_G u_gu_h(\pi(a))dh\\
&=&\int_G U_gU_h\pi(a)U_h^*U_g^*dh\\
&=&\int_GU_{gh}\pi(a)U_{gh}^*dh\\
&=&\int_GU_{h}\pi(a)U_{h}^*dh\\
&=&\int_Gu_h(\pi(a))dh=\tilde E_G(\pi(a)).
\end{eqnarray*}
Therefore, $\tilde E_G$ is an Umegaki conditional expectation with a range
\begin{equation}\label{eq:range_of_fixture}
\mathrm{Fix}(u_G):=\{x\in \pi(\mathcal A):u_g(x)=x,\,\,\forall\,g\in G\}''=\{\tilde E_G(\pi(a)):\,a\in \mathcal A\}''.
\end{equation}
We refer to \cite[Example 1.1]{St-1} for a similar exposition.
\end{rem}
\begin{lem}\label{lem:projection_Umegaki}
It holds that for all $a\in \mathcal A$,
\begin{eqnarray}
&&P_G\pi(a)P_G=\tilde E_G(\pi(a))P_G=P_G\tilde E_G(\pi(a))=P_G\tilde E_G(\pi(a))P_G,\label{eq:block1}\\
&&P_G^\perp \tilde E_G(\pi(a))P_G=P_G\tilde E_G(\pi(a))P_G^\perp=0\label{eq:block2}
\end{eqnarray}
\end{lem}
\Proof
By \eqref{eq:projection} and \eqref{eq:unitarity_for_quasi_invariant_state},
\begin{eqnarray*}
P_G\pi(a)P_G&=& \int_G\int_G U_{g_1}\pi(a)U_{g_2}dg_2dg_1  \\
&=&\int_G\int_G U_{g_1}\pi(a)U_{g_1}^*U_{g_1}U_{g_2}dg_2dg_1  \\
&=&\int_G\int_G \pi(g_1(a))U_{g_1g_2}dg_2dg_1  \\
&=&\int_G \pi(g_1(a))P_Gdg_1 \\
&=&\pi(E_G(a))P_G=\tilde E_G(\pi(a))P_G.
\end{eqnarray*}
Also,
\begin{eqnarray*}
P_G\pi(a)P_G&=& \int_G\int_G U_{g_1}\pi(a)U_{g_2}dg_2dg_1  \\
&=&\int_G\int_G U_{g_1}U_{g_2}U_{g_2}^*\pi(a) U_{g_2}dg_2dg_1  \\
&=&\int_G\int_G U_{g_1g_2}\pi(g_2^{-1}(a))dg_2dg_1  \\
&=&\int_G\int_G U_{g_1g_2}\pi(g_2^{-1}(a))dg_1dg_2  \\
&=&\int_G P_G\pi(g_2^{-1}(a))dg_2 \\
&=&P_G\tilde E_G(\pi(a)).
\end{eqnarray*}
The last equality in \eqref{eq:block1} and the equations in \eqref{eq:block2} follow by multiplying $P_G$ or $P_G^\perp$ to the left and right of the above equations.
\EndProof
\begin{thm}\label{thm:abelianness_compact_case}
Suppose that $G$ is a compact group and let $\varphi$ be a $G$-strongly quasi invariant state on a von Neumann algebra $\mathcal A$. Then, $P_G\mathcal RP_G$ is abelian if and only if $P_G\tilde E_G(\pi(\mathcal A))P_G$ is abelian, or
$P_G\mathrm{Fix}(u_G)P_G$ is abelian.
\end{thm}
\Proof It follows from Lemma \ref{lem:projection_Umegaki}  and \eqref{eq:range_of_fixture}.
\EndProof\\
Let us now come to the case that $G$ is locally compact. Suppose that $\{G_N\}_N$ is an increasing sequence of compact subgroups of $G$ such that $\cup_NG_N=G$. For each $N$, applying Lemma \ref{lem:projection_Umegaki}, we get
\begin{equation}\label{eq:projection_Umegaki_N}
P_N\pi(a)P_N=P_N\tilde E_N(\pi(a))P_N=\tilde E_N(\pi(a))P_N=P_N\tilde E_N(\pi(a)),\quad a\in \mathcal A,
\end{equation}
where we have simplified the notations as $P_N:=P_{G_N}$ and $\tilde E_N:=\tilde E_{G_N}$. By \eqref{eq:projection_range}, for each $N$ the range of $P_N$ is the closure of the vectors $\psi\in \mathcal H_\varphi$ such that $U_g\psi=\psi$ for all $g\in G_N$. Therefore, it is obvious that $\{P_N\}_N$ is a decreasing sequence of orthogonal projections, and hence it converges strongly to a limit, say $P_G$, whose range is the intersection of the ranges of all $P_N$, which is now the $U_G$-invariant set:
\[
\{\psi\in \mathcal H_\varphi:U_g(\psi)=\psi,\,\,\,\forall\,g \in G\}.
\]
\begin{lem}\label{lem:martingale_convergence}
For each $a\in \mathcal A$, the sequence $(\tilde E_N(\pi(a)))_N$ converges strongly on the range of $P_G$, that is
\begin{equation}\label{eq:strong_limit}
 \lim_{N\to \infty} \tilde E_N(\pi(a))P_G\psi=P_G\pi(a)P_G\psi, \quad \psi\in \mathcal H_\varphi.
\end{equation}
\end{lem}
\Proof
For any $\psi\in \mathcal H_\varphi$, by \eqref{eq:projection_Umegaki_N}, we have
\begin{eqnarray*}
 \tilde E_N(\pi(a))P_G\psi&=& \tilde E_N(\pi(a))P_N(P_G\psi)\\
 &=&P_N\pi(a)P_G\psi\\
 &\to&P_G\pi(a)P_G\psi\quad \text{as }N\to \infty
 \end{eqnarray*}
 since $P_N\to P_G$ strongly.
 \EndProof\\
\begin{rem}\label{rem:Umegaki_locally_compact}
We denote the strong limit of $(\tilde E_N(\pi(a))P_G)_N$ by $\tilde E_G(\pi(a))P_G$. By \eqref{lem:martingale_convergence}, it is equal to $P_G\tilde E_G(\pi(a))P_G$. One can show that  $P_G\tilde E_G(\cdot)P_G$ is an Umegaki conditional expectation onto $P_G\mathrm{Fix}(u_G)P_G$, where
\[
\mathrm{Fix}(u_G):=\{x\in \pi(\mathcal A):u_g(x)=x,\,\,\forall\,g\in G\}''.
\]
\end{rem}
 \begin{thm}\label{thm:projection_Umegaki_locally_compact_group}
 Let $G$ be a locally compact automorphism group acting on a von Neumann algebra $\mathcal A$, for which there is an increasing sequence of compact subgroups $\{G_N\}_N$ such that $\cup_NG_N=G$. Let $\varphi$ be a faithful and normal $G$-strongly quasi invariant state on $\mathcal A$. With the notations given above we have that $P_G\mathcal RP_G$ is abelian if and only if $P_G\mathrm{Fix}(u_G)P_G$ is abelian.
\end{thm}
\Proof It follows from Lemma \ref{lem:martingale_convergence} and Remark \ref{rem:Umegaki_locally_compact}.
\EndProof
\begin{rem} Notice that if $\mathrm{Fix}(u_G)$ is $\mathbb C \unit$, then $P_G\mathrm{Fix}(u_G)P_G$ is one dimensional and hence it is automatically abelian.
\end{rem}
\noindent{\bf Acknowledgement}. A. Dhahri is a member of GNAMPA-INdAM and he has been supported by the MUR grant Dipartimento di Eccellenza 2023-2027 of Dipartimento di Matematica, Politecnico di Milano. The work of H. J. Yoo was supported by the National Research Foundation of Korea (NRF) grant funded by the Korean government (MSIT) (No. RS-2023-00244129).

 \end{document}